\documentclass[11pt]{amsart}
\usepackage{graphicx}
\usepackage{amssymb}
\def\L{{\Lambda}}

\newtheorem{theorem}{Theorem}
\newtheorem{definition}{Definition}

\title[Quasiconformal maps with H\"older continuous dilatation]{A note on quasiconformal maps with H\"older-continuous dilatation}
\author{James T. Gill}
\address{Department of Mathematics and Computer Science\\ Saint Louis University\\ St. Louis, MO 63103}
\email{jgill5@slu.edu}
\author{Steffen Rohde}
\address{Department of Mathematics\\ University of Washington\\ Seattle, WA 98195-4350}
\email{rohde@math.washington.edu}

\thanks{The authors were partially supported by the following NSF grants: the first author, \#1004721; the second author, \#1068105.}
\subjclass[2010]{Primary 30C65}
\keywords{bi-Lipschitz, quasiconformal, H\"older continuous, conformal modulus}

\begin{document}

\begin{abstract}
Recently J. Mateu, J. Orobitg, and J. Verdera showed that a H\"older continuous complex dilatation supported on smooth domains is a sufficient condition for the resulting quasiconformal map to be bi-Lipschitz.  Their proof is analytic and based on properties of the Beurling-Ahlfors transform.  We give an alternate, more geometric proof and use it to extend their result to supporting domains with positive angle corners.
\end{abstract}

\maketitle

\section{Introduction}

The Beltrami equation
\[ \frac{\partial \Phi}{\partial \overline{z}} = \mu(z) \frac{\partial \Phi}{\partial z}, \hspace{.5cm}\mbox{a.e. } z \in \mathbb{C} \]
where $\mu$ is Lebesgue measurable in $\mathbb{C}$, is an elliptic (uniformly elliptic if $\| \mu \|_\infty<1$) partial differential equation whose homeomorphic solutions are the quasiconformal maps (finite distortion maps if $\| \mu \|_\infty =1$) of the plane.  The quasiconformal maps enjoy a variety of important properties (see  the recent book \cite{AIM} for a modern treatment from the point of view of PDE theory or the classic treatise \cite{LV}) and are often a useful tool for situations where conformal maps turn out to be too rigid.  They are, however, not a perfect stand-in.  They fail to be locally bi-Lipschitz, for example, and this is the topic of our note.

A necessary and sufficent condition for a quasiconformal map to be bi-Lipschitz is unknown.  Perhaps the first step of understanding is due to J. Schauder who showed that it is sufficient that $\mu$ be H\"older continuous and compactly supported (see Chapter 15 of \cite{AIM}).  In \cite{MOV}, J. Mateu, J. Orobitg, and J. Verdera prove
\begin{theorem}
Let $\{ \Omega_j \}$, $1 \leq j \leq N$, be a finite family of disjoint bounded domains of the plane with boundary of class $C^{1 + \varepsilon}$, $0 < \varepsilon < 1$, and let $\mu = \sum_{j=1}^n \mu_j \chi_{\Omega_j}$, where $\mu_j$ are $\varepsilon$-H\"older continuous functions on $\Omega_j$ and $\| \mu_j \|_\infty < 1$ for each $j = 1, \ldots N$.  Then the  principal solution associated to the $\mu$-Beltrami equation (i.e. the one which is $z+O(1/z)$ near $\infty$) is a bi-Lipschitz quasiconformal map.
\end{theorem}
By $\varepsilon$-H\"older continuous on $\Omega$ we mean that there exists a $C>0$ so that 
\[ |\mu(x) - \mu(y) | \leq C |x-y|^\varepsilon \]
for all pairs $x,y \in \Omega$.

Informally, Theorem 1 extends the result of Shauder to allow the dilatation to have jumps, provided that the jumps occur at reasonably smooth curves.  In this note we give an alternate proof of this theorem using the more geometric methods of \cite{R}.  In addition we will also show that our geometric methods extend Theorem 1 to domains with (properly defined) corners if the dilatation takes certain allowable values (which depend on the angle at the corner).

\begin{theorem}
The domains in Theorem 1 may be taken to each have finitely many corners, provided the dilatation takes certain allowable values at each corner.
\end{theorem}

We deliberately postpone the definition of both the allowable values of $\mu$ and the type of corners allowed until later, when they will hopefully appear natural to the reader.

The novelty of this theorem is that even the most trivial examples of domains with corners prove to be quite complicated.  For instance in \cite{AC}, A. Ch\'eritat considers the case of the Beltrami equation which is a constant on the unit square and zero elsewhere.  The computation of these maps and the image of the unit square under them are non-trivial.  This is in direct contrast with the case of a Beltrami equation supported on the unit disk, as discussed in Sections 2.4 and 2.5 below.

The structure of this note is as follows: in Section 2 we outline some examples and tools which will be useful in our proofs, in Section 3 we will give an alternate proof of the theorem of Mateu, Orobitg, and Verdera, and in Section 4 we extend this result to domains with corners.

\subsection{Acknowledgement}
Part of this work was done while the second author was visiting the  Centre de Recerca Matem\`atica, and he would like to thank the CRM for their support and their hospitality. He would also like to  thank Tadeusz Iwaniec for stimulating discussions related to this work.  While this research was done the first author was supported by the University of Washington Department of Mathematics.

\section{Tools and Examples}
In this section we describe some standard properties of conformal moduli which we will use in our proofs and discuss the role of both the Koebe distortion theorem and the Lehto integral condition in what follows.  Finally, we describe specific examples of bi-Lipschitz quasiconformal maps which we will use in our proofs.  The first example is rather simple, the second is more delicate.

\subsection{Modulus Estimates}

Modulus is a conformal invariant for doubly connected domains in the plane.  Here we provide a brief overview of this invariant and discuss some of the properties which will be of use to us.  Suppose we have an annulus centered at 0: $A(0,r_1,r_2): = \{ z: 0 \leq r_1 < |z| < r_2 \leq \infty \}$.  The positive quantity
\[ M(A) = \log \frac{r_2}{r_1} \]
is called the  \textit{modulus} of $A$.  For other doubly connected domains $B$, there is a (round) annulus centered at zero which is the conformal image of $B$.  We define the \textit{modulus} of $B$, $M(B)$, to be the modulus of such an annulus.  It can be shown that $M(B)$ is well defined.  This notion gives a simply stated geometric definition for quasiconformal maps, which can be shown (see Chapter I of \cite{LV}) to be equivalent to the other common definitions.

\begin{definition}
A {\em K-quasiconformal map} $f$ is a sense-preserving homeomorphism of a domain $D$ with the property that all doubly connected $B$ with $\overline{B} \subset D$ and their images, $f(B)$, satisfy
\[ M(f(B)) \leq K \cdot M(B). \]
\end{definition} 

The calculation of the modulus of an arbitrary doubly connected region is often difficult and the image of an annulus under a quasiconformal map may no longer be an annulus.    We will get around this difficulty by showing here that if an annulus has modulus large enough, then its image under a global $K$-quasiconformal map is quantitatively ``almost round''.

Let $f$ be a $K$-quasiconformal homeomorphism of $\mathbb{C}$.  Such global $K$-quasiconformal maps are \textit{quasisymmetric}, meaning there exists an increasing homeomorphism $\eta :[0, \infty) \to [0, \infty)$ such that for all $z_0, z_1, z_2 \in \mathbb{C}$,
\[ \frac{ |f(z_0)-f(z_1)|}{|f(z_0)-f(z_2)|} \leq \eta \left( \frac{|z_0-z_1|}{|z_0-z_2|} \right) . \]
In fact, $\eta=\eta_K$ depends only on $K$ (see Chapter 3 of \cite{AIM}).

\begin{figure}[htp]
\centering
\includegraphics[scale=0.50]{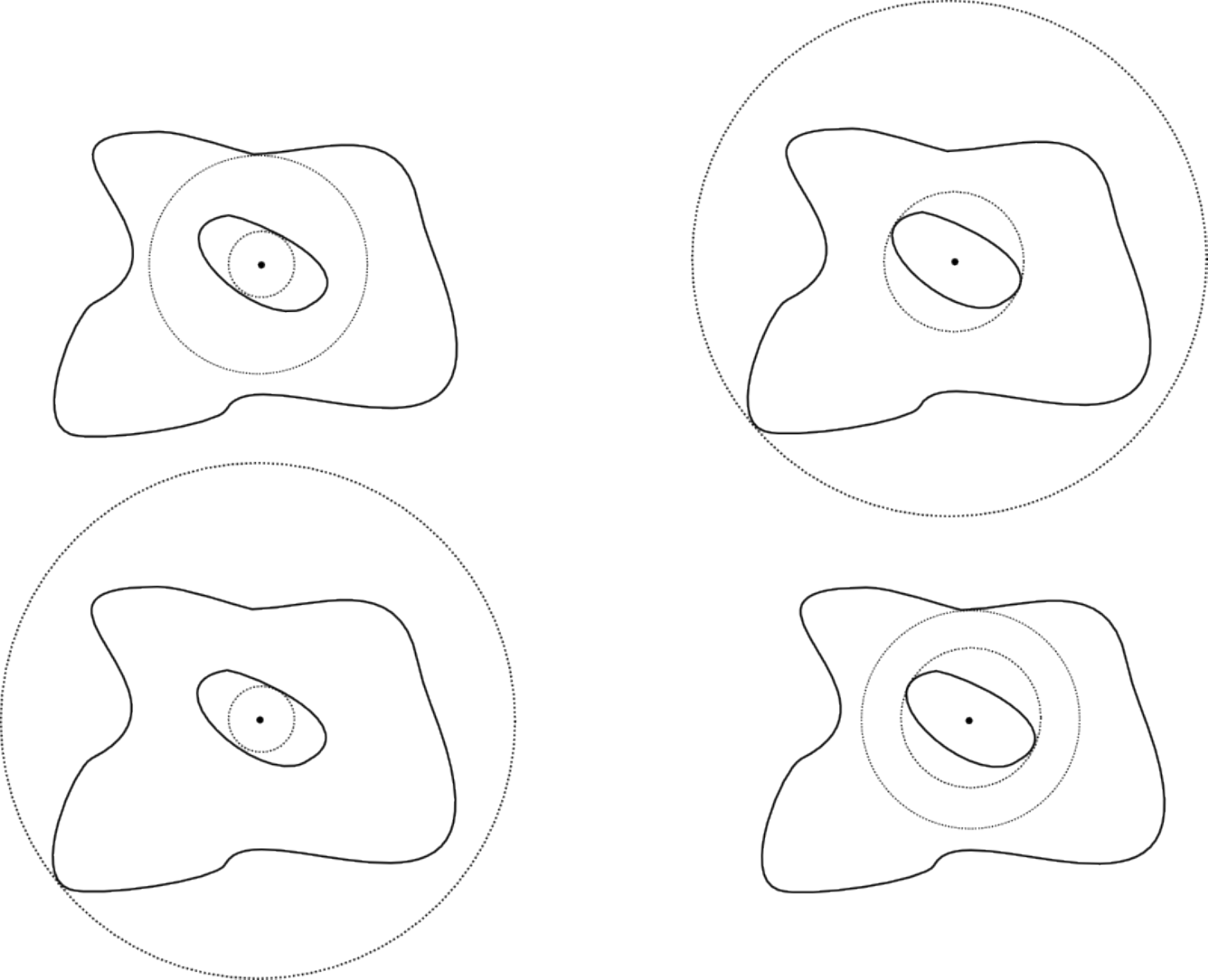}
\caption{A sketch of annuli $B,C,D,$ and $E$, clockwise from the upper left, around a single region $f(A)$}\label{fig:annuli}
\end{figure}

Consider the annulus $A = A(z_0,r_1,r_2)$ with $r_1 < r_2$  and $\eta(r_1/r_2) < 1$ (this is the ``large enough'' condition).  Let
\[ \rho_1 := \min_{|z-z_0|=r_1}|f(z)-f(z_0)| \mbox{ and } \rho_2 := \min_{|z-z_0|=r_2} |f(z)-f(z_0)| \]
and $R_1$ and $R_2$ given by replacing $\min$ with $\max$ above, respectively.
As $\eta(r_1/r_2) < 1$, we have $R_1 < \rho_2$.  Let us define the annuli
\[ B=A(f(z_0),\rho_1, \rho_2) ,\hspace{1cm}  C=A(f(z_0), R_1, R_2), \]
\[ D=A(f(z_0),\rho_1, R_2), \hspace{1cm} E=A(f(z_0),R_1, \rho_2). \hspace{.1cm} \]
A sketch of such annuli around a doubly connected region representing $f(A)$ is shown in Figure \ref{fig:annuli}.

Note that $E$ is the annulus of largest modulus centered at $f(z_0)$ contained in $f(A)$ and $D$ is the annulus of smallest modulus centered at $f(z_0)$ which contains $f(A)$.  Using the definition of quasisymmetry, the monotonicity of the modulus with respect to inclusion, and the fact that $r_2/r_1$ is large we get the following chain of inequalities:
\begin{equation*}
0 \leq \log \left( \frac{1}{\eta_K (r_1/r_2)} \right) - \log \eta_K (1) \leq M(C) - \log \eta_K(1) \leq \log \frac{\rho_2}{R_1} = M(E) \hspace{2cm}
\end{equation*}
\begin{equation}
\hspace{2cm} \leq M(f(A)) \leq M(D) = \log \frac{R_2}{\rho_1} \leq  \log \eta_K \left( \frac{r_2}{r_1} \right). \label{modineq}
\end{equation}
Of particular importance will be the difference between $M(D)$ and $M(E)$:
\begin{equation}
 0 \leq M(D)-M(E) = \log (R_2/\rho_1) - \log (\rho_2/R_1) \leq 2 \log (\eta_K(1)) :=C_K \label{annuliineq}.
\end{equation}
This will be especially useful in concert with the following theorem due to the second author in \cite{R} which provides a geometric criterion for quasiconformal bi-Lipschitz maps.

\vspace{2mm}\noindent{\bf Theorem A.}
\textit{
Let $f$ be a $K$-quasiconformal homeomorphism of $\mathbb{C}$ and $E \subset \mathbb{C}$ any set.  Assume that there is a constant $N$ such that the difference in conformal modulus of $A$ and $f(A)$ is bounded by $N$, i.e.
\begin{equation} | M(f(A)) - M(A) | \leq N  \label{gbl} \end{equation}
holds for all annuli $A$ centered at points of $E$ with the property that both boundary circles meet $E$.  Then the restriction of $f$ to $E$ is bi-Lipschitz.}
\vspace{1mm}

When we apply this result we will show this condition only for annuli with large modulus, as the definition of quasiconformal maps above already gives (\ref{gbl}) for annuli with small modulus.

\subsection{The Koebe Distortion Theorem}
The Koebe distortion theorem (see Chapter 2 of \cite{AIM}) is a powerful tool from conformal mapping theory.  It implies that if $\Lambda$ is a compact set contained in a domain $\Omega$, then there is an $M>1$, dependent on $\Lambda$ inside $\Omega$, such that for any conformal map $f$ of $\Omega$ and any pair of points $z,w \in \Lambda$,
\[ \frac{1}{M} |f'(w)| \leq |f'(z)| \leq M |f'(w)|. \]
Hence the modulus of the derivative of $f$ at $z$ cannot be far from the modulus of the derivative at $w$, and from this we can see that conformal maps are locally bi-Lipschitz.  Now consider two quasiconformal maps $f_1,f_2$ of a common domain $\Omega$ with $\mu_{f_1} = \mu_{f_2}$ on $\Omega$.  Assume that $f_2$ is bi-Lipschitz on a compact set $\Lambda \subset \Omega$.  The composition formula for dilatations (which we will use often in the sequel) says that
\begin{equation} \mu_{f_1 \circ f_2^{-1}} (\zeta) = \frac{ \mu_{f_1} (z) - \mu_{f_2}(z)}{ 1 - \mu_{f_1}(z) \overline{\mu_{f_2} (z) } } e^{2 i \arg \partial f_2 (z)}  \label{mucomp} \end{equation}
where $\zeta = f_2(z)$.  Thus $f_1 \circ f_2^{-1}$ is bi-Lipschitz on $f_2(\Lambda)$ since it is conformal,  and $f_1 = (f_1 \circ f_2^{-1}) \circ f_2$ is bi-Lipschitz on $\Lambda$ as the composition of two bi-Lipschitz maps.  So maps with identical dilatation on a domain share the local bi-Lipschitz property.  This will allow us to compare known bi-Lipschitz solutions to the Beltrami equation with arbitrary solutions.

\subsection{The Lehto Inequality}

Another important tool for us will be an integral inequality originally discovered by O. Lehto in \cite{L} but is perhaps more easily found in Chapter V of \cite{LV}.   Let $f$ be a $K$-quasiconformal map on some domain containing the annulus $A=A(x,r,R)$.  Let $\mu_f$ be the complex dilatation of $f$. Then
\begin{equation}
|M(f(A)) - M(A)| \leq C(K) \int_A \frac{|\mu_f(y)|}{|x-y|^2}  dy, \label{lehto}
\end{equation}
where $dy$ here stands for 2-dimensional Lebesgue measure. Considering the criterion described in Theorem A above, it is clear that (\ref{lehto}) may prove to be a useful tool.  In fact, the singularity of the integral in (\ref{lehto}) already hints at the importance of H\"older continuity of $\mu$ in what follows.

\subsection{A Simple Bi-Lipschitz Example}

In order to apply the Lehto Inequality (\ref{lehto}) for all annuli around a given point, one needs $\mu$ to vanish at that point in such a way to eliminate the singularity in the integral of (\ref{lehto}).  This will be done through composing with known bi-Lipschitz maps.  We start with a simple known example where $\mu$ is compactly supported. Let $c \in \mathbb{D}$ be given.  Consider the map
\begin{equation} f_c(z) := \left\{ \begin{array}{cl}
z + \frac{c}{z} & \mbox{ in } \mathbb{D}^c\\
z + c \overline{z} & \mbox{ in } \overline{\mathbb{D}} \end{array} \right. \label{simple} \end{equation}
This map is a homeomorphism of the plane and it can easily be seen to be bi-Lipschitz.  Note that the image of the unit circle is an ellipse and the map is conformal outside the unit disk.  It is also the unique solution to the Beltrami equation with dilatation $c \cdot \chi_\mathbb{D}$ and $f(z) = z + O(1/z)$, often called the principal solution.

\subsection{A More Complicated Bi-Lipschitz Example}

While the solution to the constant dilatation Beltrami equation supported on the disk is always bi-Lipschitz, and this is also true of the half-plane, this is radically not true of sectors of angle $\not=  \pi$.  We show this by constructing such a map.

Choosing the pricipal branch of the logarithm and $\theta_0 \in [0,\pi) \cup (\pi, 2\pi)$, we let
\[ \mathbb{S}_0^{\theta_0} = \{z \not= 0 : 0 \leq \arg z <\theta_0 \}, \hspace{5mm}\mathbb{S}_{\theta_0}^{2\pi} = \{z \not= 0 : \theta_0 \leq \arg z <2\pi \}. \] 
We wish to find a solution to the Beltrami equation where $\mu$ is given by
\[ \mu = \left\{ \begin{array}{cl}
c & \mbox{ in } \mathbb{S}_0^{\theta_0}\\
0 & \mbox{ in } \mathbb{S}_{\theta_0}^{2\pi} \end{array} \right. . \]
We start by setting
\[ f_1(z) = \frac{z + c \overline{z}}{1+c} \mbox{ in } \mathbb{S}_0^{\theta_0}.\]
This conveniently gives $f_1(x) = x$ for real numbers $x>0$ and $\mu_{f_1}=c$ in $\mathbb{S}_0^{\theta_0}$.  To extend $f_1$ continuously to the second sector we examine
\[ f_1(r e^{i \theta_0}) = \frac{e^{i \theta_0} + c e^{- i \theta_0}}{1+c} r. \]
Let $Re^{i \theta_1} := \frac{e^{i \theta_0} + c e^{- i \theta_0}}{1+c}$ for $R>0$ and $0 \leq \theta_1 < 2\pi$, and $a := \frac{Re^{i \theta_1}}{e^{i \theta_0}}$.  Now we set
\[ f_1(z) = a z \mbox{ in } \mathbb{S}_{\theta_0}^{2\pi}.\]
This $f_1$ would be our solution, but the values may not match up on the two sides of the positive real axis.  Our solution to this will be to slice open the plane via a logarithm, then apply a rotation and dilation to get these two parts to match up continuously, and then we sew the plane back together with an exponential.

\begin{figure}[htp]
\centering
\includegraphics[scale=0.40]{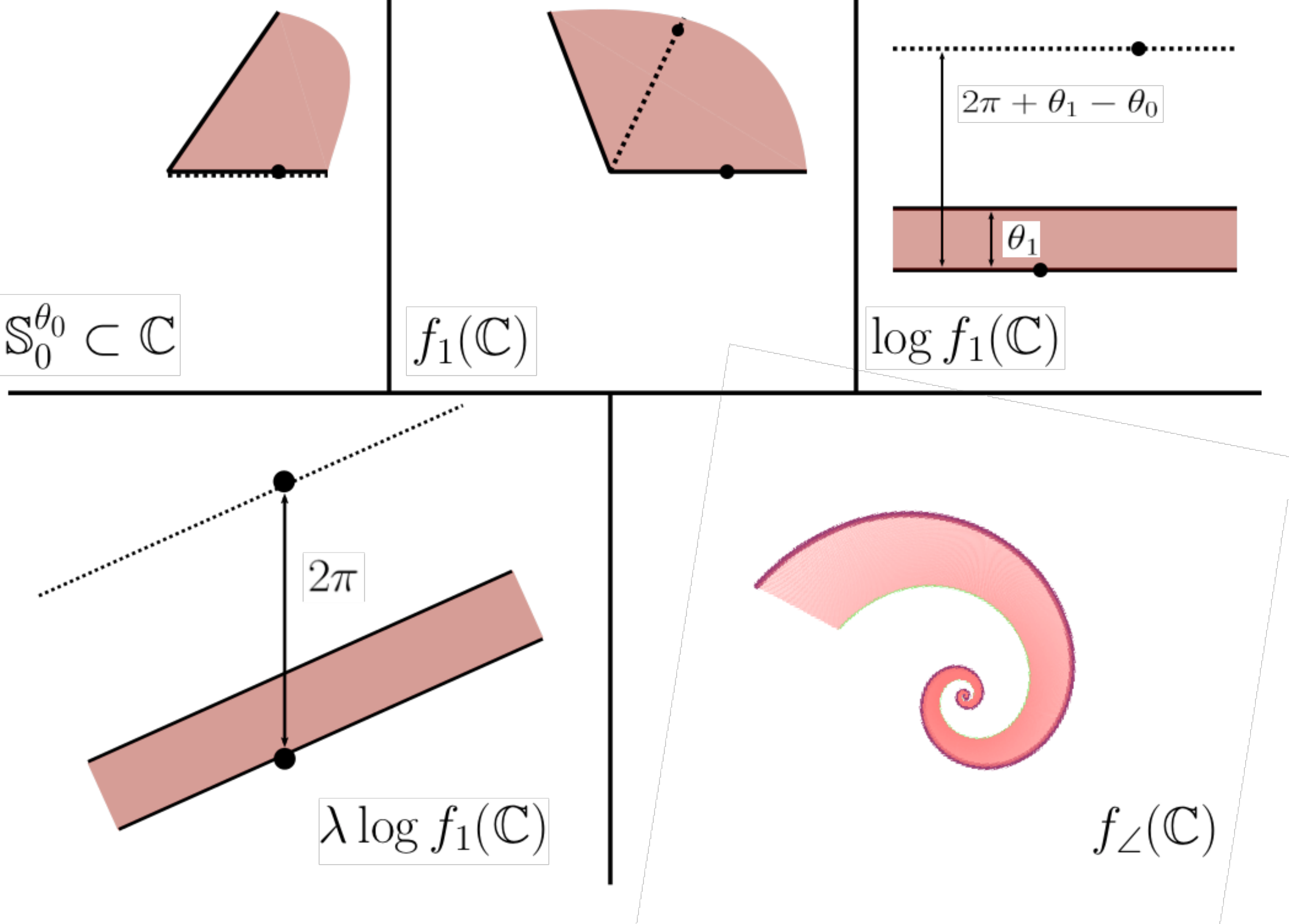}
\caption{The components of the map $f_\angle$}\label{fig:test2}
\end{figure}

Using the principal branch of the logarithm we consider the image of the two sides of $\mathbb{R}^+$:
\begin{eqnarray*}
\log (f_1 (r)) & = & \log r\\
\log (f_1 (r e^{i 2\pi})) & = & \log r + \log R + i(2\pi + \theta_1 - \theta_0)
\end{eqnarray*}
In order to line up these two sides of the positive real line, we multiply by
\[ \lambda = \lambda(c, \theta_0) := \frac{2 \pi i}{\log R + i(2\pi + \theta_1 - \theta_0)} \]
and our Beltrami solution is
\[ f_\angle (z) = e^{\lambda \log (f_1 (z))}, \]
so that the image of the sector $S_0^{\theta_0}$ is a logarithmic spiral. See Figure \ref{fig:test2} for a visual description of this map. 

This map will be bi-Lipschitz if and only $Re (\lambda) = 1$.  As $\lambda$ depends on $c$ and $\theta_0$, it follows that when an angle for the sector is fixed, only certain values for the dilatation in that sector give rise to bi-Lipschitz maps.  

As maps with identical dilatation on a domain are simultaneously locally bi-Lipschitz by the discussion in Section 2.2, any description of sufficient conditions for a quasiconformal map to be bi-Lipschitz  must take the rigidity of this example into account. 
\begin{definition}Given an angle $\theta_0 \in [0,2\pi)$, we say that $c \in \mathbb{D}$ is {\em allowable for $S_0^{\theta_0}$} if $Re(\lambda(\theta_0,c)) = 1$.
\end{definition}
The reader interested in dilatation supported on the unit square would find the well-illustrated note \cite{AC} interesting.


\section{Geometric Proof of the Theorem of Mateu, Orobitg, and Verdera}

Our proof of Theorem 1 begins with the simplifying assumption that $j = 1$  (see Section 7 of \cite{MOV} for the standard details). Let $\Omega$ be a bounded domain in $\mathbb{C}$ with $\partial \Omega$ a $C^{1+\varepsilon}$ curve and $\mu$ be $\varepsilon$-H\"older continuous in $\Omega$ and $\mu \equiv 0$ in $\Omega^c$.  Let $\Phi$ be the solution of $\overline{\partial} \Phi = \mu \partial \Phi$ with $\Phi(z)= z +O(1/z)$ near $\infty$.

\subsection{Reduction to the Unit Disk}
We wish to reduce Theorem 1 to the case where $\Omega = \mathbb{D}$, the unit disc.  Let $\phi$ be a conformal map from $\overline{\mathbb{D}}^c$ to $\overline{\Omega}^c$ with $\phi(z) = z + O(1)$ near $\infty$ and let $\psi$ be a conformal map from $\mathbb{D}$ to $\Omega$.  By the Kellog-Warschawski theorem on the boundary correspondence of conformal maps (see Chapter 3 of \cite{P}), both maps extend to $\partial \mathbb{D}$ as $C^{1+\varepsilon}$ mappings with non-vanishing derivatives. So both maps extend to $\partial \mathbb{D}$ as bi-Lipschitz conformal maps.  Then $\psi^{-1} \circ \phi : \partial \mathbb{D} \to \partial \mathbb{D}$ is a $C^{1+\varepsilon}$ homeomorphism with a derivative which is never 0.  By a pre-rotation of $\psi$, assume that $\psi^{-1} \circ \phi (1) = 1$.  We wish to extend $\psi^{-1} \circ \phi$ to the whole unit disc so that the extension is bi-Lipschitz, quasiconformal, and with $\varepsilon$-H\"older continuous dilatation $\mu$.  Set $\gamma(\theta) = \arg{(\psi^{-1}(\phi(e^{i \theta})))}$ where the argument takes values in $[0,2\pi)$.
A simple  bi-Lipschitz extension would be given by $re^{i\theta}\mapsto  r e^{i \gamma(\theta)},$ but its Beltrami coefficient would not even be continuous at $0.$ For this reason  we define our extension as follows:
\[ F(re^{i\theta}) = r e^{i\left[ r \gamma(\theta) + (1-r) \theta \right]}. \]
Clearly, $F$ is a homeomorphism of $\overline{\mathbb{D}}$.
Using polar coordinates we calculate
\[ \mu_F (r e^{i \theta}) = \frac{r(1-\gamma'(\theta)) + i r (\gamma(\theta) - \theta)}{2 - r(1-\gamma'(\theta)) + i r(\gamma(\theta)-\theta)}\  e^{2i\theta}\ . \]
From these formulae it is not hard to check that $F$ is indeed bi-Lipschitz, quasiconformal, and has $\varepsilon$-H\"older continuous dilatation.   Let
\[ f = \left\{ \begin{array}{cr}
\phi & \mbox{ outside } \mathbb{D}\\
\psi \circ F & \mbox{ inside } \mathbb{D}
\end{array} \right. .\]
Then $f$ is a bi-Lipschitz map from $\mathbb{C}$ to itself with $\varepsilon$-H\"older continuous $\mu$ supported on the unit disk, and the image of the unit disk under $f$ is $\Omega$. This $f$ allows us to reduce Theorem 1 to the case of $\Omega = \mathbb{D}$, as $\Phi \circ f$ fulfills the assumptions of Theorem 1 with $\Omega = \mathbb{D}$ (note that via the composition rule (\ref{mucomp}) the criterion on $\varepsilon$-H\"older continuity of $\mu_{\Phi \circ f}$  still holds), and because if $\Phi \circ f$ is bi-Lipschitz, then so is $\Phi$.

\subsection{Proof on the Unit Disk}

Let $f$ be a homeomorphism satisfying $\overline{\partial} f = \mu_f \partial f$ almost everywhere in $\mathbb{C}$ with $\mu_f$ an $\varepsilon$-H\"older continuous function in $\mathbb{D}$, $\mu_f \equiv 0$ outside $\mathbb{D}$, 
$k=\|\mu_f\|_\infty<1$
and $f = z + O(1/z)$ near $\infty$.  The Koebe distortion theorem automatically gives that $f$ is bi-Lipschitz outside, say, $\frac32\mathbb{D}$.  We will show that $f$ is bi-Lipschitz in $\mathbb{C}$ by applying Theorem A with $E = 2 \mathbb{D}$.  First, as pointed out in Section 2.1, $K-$quasiconformal maps distort the modulus of annuli by a factor less than or equal to $K$ so we may assume that $A = A(x_0,r,R)$ with $x_0 \in 2\mathbb{D}$ has $M(A)$ large.  We must find an upper bound for $|M(f(A)) - M(A)|$ independent of the specific annulus $A$.  If $A \cap \mathbb{D} = \emptyset$, then $M(f(A)) = M(A)$ as  the modulus is a conformal invariant.

Let $A \cap \mathbb{D}$ be non-empty with $x_0 \in 2\mathbb{D}$.  In order to apply the Lehto condition (\ref{lehto}), we first reduce to the case $\mu(x_0) = 0$ by precomposing $f$ with an appropriate bi-Lipschitz $f_c$ from Section 2.4.  Set
\[ x^* := \left\{ \begin{array}{cl}
x_0 & \mbox{if } x \in \mathbb{D}\\
\frac{x_0}{|x_0|} & \mbox{else}
\end{array} \right. \]
and $\mu_0(\cdot) = \mu_f(x^*) \cdot \chi_{\mathbb{D}}(\cdot)$, where $\mu_f$ is understood to take its unique limit value on $\partial \mathbb{D}$ via approaching the boundary from the inside as $\mu$ is uniformly H\"older continuous in $\mathbb{D}$.  Let $c = \mu_f(x^*)$ and $f_c$ the bi-Lipschitz principal solution to $\overline{\partial} f_c = \mu_0 \partial f_c$ given in Section 2.4.  We now consider the map $\tilde f = f \circ f_c^{-1}$.  We wish to show $\tilde f$ is bi-Lipschitz as well.  If $y \notin f_c(\mathbb{D})$, then $\mu_{\tilde f}(y)=0$.  Otherwise, $y = f_c(x)$ for some $x \in \mathbb{D}$.  Via the composition formula for dilatations (\ref{mucomp}) and the $\varepsilon$-H\"older continuity of $\mu_f$
\[ (1-k^2)|\mu_{\tilde f} (y) | \leq |\mu_f(x) - \mu_{f_c}(x)| = |\mu_f(x) - \mu_f(x^*)| \leq C|x-x^*|^\varepsilon \leq C |x-x_0|^\varepsilon. \]
Since $f_c$ is bi-Lipschitz, $|x-x_0| \leq L |y-f_c(x_0)|$ for some $L = L(\|\mu_f\|_\infty).$  Thus
\begin{equation} |\mu_{\tilde f}(y) | \leq C' |y-f_c(x_0)|^\varepsilon. \label{mus} \end{equation}
Consider the doubly connected region $f_c(A)$ and the related annuli $E$ and $D$ defined by $f_c(A)$ via Section 2.1. Then via the triangle inequality
\[
|M(f(A)) - M(A)| \leq \hspace{7cm}
\]
\begin{equation} \label{trimod}
           \hspace{1.5cm} |M(f(A)) - M(\tilde f(E))| + |M(\tilde f(E)) - M(E)| + |M(E) - M(A)|.
\end{equation}

We work on each of these quantities on the right hand side of (\ref{trimod}) separately.  Let $k=\|\mu_f\|_\infty$ and $K = \frac{1+k}{1-k}$.  The second quantity, $|M(\tilde f(E)) - M(E)|$, can be estimated by using the Lehto inequality (\ref{lehto}) in Section 2.3:
\begin{equation} |M(\tilde f(E)) - M(E)| \leq C(K) \int_E \frac{|\mu_{\tilde f} (y)|}{|f_c(x_0) - y|^2}dy \leq N
\label{lehto2} \end{equation}
by (\ref{mus}) where the constant N depends on $\|\mu_f\|_\infty$ and the $\varepsilon$-H\"older constant of $\mu_f$, but not on $x_0$ or $E$.

The first quantity on the right hand side of (\ref{trimod}) can be rewritten
$|M(\tilde f f_c (A)) - M(\tilde f(E))|$. By monotonicity of the modulus, $M(\tilde f f_c (A)) \geq M(\tilde f(E))$, and $M(\tilde f(D)) \geq M(\tilde f f_c(A))$ so
\[|M(\tilde f f_c (A)) - M(\tilde f(E))| \leq \hspace{7cm}
\]
\[ \hspace{1.5cm} |M(\tilde f(D)) -M(D)| + |M(D)-M(E)| + |M(E)-M(\tilde f(E))|. \]
The first and last terms are bounded by a constant from the Lehto inequality (\ref{lehto}) as above in (\ref{lehto2}).  By (\ref{annuliineq}) the middle term is bounded by a constant as well. So $|M(\tilde f f_c (A)) - M(\tilde f(E))| \leq 2N+ C_K$.

For the third quantity in (\ref{trimod}), $|M(E) - M(A)|$, we again use the triangle inequality to get
\[ |M(E)-M(A)| \leq |M(E)-M(f_c(A))| + |M(f_c(A))-M(A)|. \]
The first is less than a constant depending on $K$ by (\ref{annuliineq}) and the second is bounded by a constant depending on $K$ by using the geometric criterion for bi-Lipschitz maps from \cite{R}, as $f_c$ is bi-Lipschitz.

Taken altogether, $|M(A) - M(f(A))| \leq P$, for some constant $P$ depending only on $\|\mu_f\|_\infty$ and the $\varepsilon$-H\"older constant of $\mu_f$. So by Theorem A, a quasiconformal map with $\varepsilon$-H\"older continuous $\mu$ supported on $\mathbb{D}$ is bi-Lipschitz.  We reduced Theorem 1 to this case, so Theorem 1 is proved.

\section{Domains with Corners}

\subsection{Definitions and Precise Statement of Theorem 2}One of the nice outcomes of this geometric proof to Theorem 1 is that its methods are amenable to domains without smooth boundary, as long as we have some knowledge of the behavior of conformal maps on such domains. The theory of O. Kellogg and S. Warschawski discussed in Chapter 3 of \cite{P} will give us the needed tools.  We state the definition of a corner from \cite{P}:

\begin{definition} Let $\Omega$ be any bounded simply connected domain in $\mathbb{C}$ with a locally connected boundary and $f:\mathbb{D} \to \Omega$ be a conformal map onto $\Omega$ which is continuous in $\overline{\mathbb{D}}$.  Let $\zeta = e^{i \theta} \in \partial \mathbb{D}$. We say that $\partial \Omega$ has a {\em corner of opening} $\alpha$, $0 \leq \alpha \leq 2$, at $f(\zeta)$ if
\[ \arg [f(e^{i t}) - f(e^{i \theta})] \to  \left\{ \begin{array}{cl}
\beta & \mbox{ as } t \to \theta+,\\
\beta + \alpha \pi & \mbox{ as }t \to \theta-
\end{array} \right. \]
for some $\beta$, $0 \leq \beta \leq 2\pi$. In addition, we say that $\partial \Omega$ has a {\em H\"older-continuous corner} at $f(\zeta)$ if there are closed arcs $A^{\pm} \subset \partial \mathbb{D}$ ending at $\zeta$ and lying on opposite sides of $\zeta$ that are mapped onto $C^{1+\varepsilon}$ Jordan arcs $C^+$ and $C^{-}$ forming a corner of opening $\alpha$ at $f(\zeta)$ for some $\varepsilon > 0$.
\end{definition}

If the boundary of a domain $\Omega$ is composed of finitely many $C^{1+\varepsilon}$ closed arcs meeting at corners, then each corner is H\"older-continous.  Domains with H\"older-continuous corners are nice because of the following theorem due originally to Warschawski but found as Exercise 3.4.1 in \cite{P}:

\vspace{2mm}\noindent{\bf Theorem B.}
\textit{
Let $\partial \Omega$ have a H\"older continous corner of opening $\alpha > 0$ at $f(\zeta)$ with H\"older constant $\varepsilon >0$, where $f$ is a conformal map from $\mathbb{D}$ onto $\Omega$.  Then, for some $a \not=0$,
\begin{equation} f(z) = f(\zeta) + a(z-\zeta)^{\alpha} + O(|z-\zeta|^{(\alpha + \alpha \varepsilon)}) \label{cornerexp}
\end{equation}
as $z \to \zeta$ in $\overline{\mathbb{D}}$.
}
\vspace{2mm}

We are almost ready to state and prove our theorem concerning domains with corners, but we must first be careful of the fact that pre-composing a quasiconformal map with a conformal map does not leave the dilatation invariant.  We refine our definition of allowable dilatations taking into account the composition formula for complex dilatations (\ref{mucomp}):

\begin{definition}
Let $\Omega$ be a bounded domain with a corner of opening $ \alpha$ at $z_0 = f(\zeta)$.  Let $f(z) = e^{-i  \beta}(z-z_0)$ be the conformal map of the plane that moves the corner to 0 and one of the arcs forming the corner to an arc tangent with the real line at 0 and the other tangent to a ray at positive angle $\alpha \pi$ from the first.  We say that $c \in \mathbb{D}$ is {\em allowable for $(\Omega, z_0)$} if $c \cdot e^{-2i \beta}$ is allowable for $\mathbb{S}_0^{\alpha \pi}$.
\end{definition}

We are now ready for the precise statement of our theorem concerning domains with corners:

\vspace{2mm}\noindent{\bf Theorem 2.}
\textit{
Let $\Omega$ be a bounded domain whose boundary consists of $N$ closed sub-arcs of $C^{1+\varepsilon}$ arcs.  Let the endpoints of these sub-arcs be labeled $z_j$ and let $\partial \Omega$ have a corner of opening $\alpha_j \in (0,2)$ at $z_j$.  Assume the $z_j$ are distinct.  Let $\mu$ be an $\varepsilon$-H\"older continuous function on $\Omega$ with the values $\mu(z_j)$ each allowable for $(\Omega, z_j)$ and $\mu \equiv 0$ outside $\Omega$.  Assume, in addition, that $\| \mu \|_\infty < 1$.  Then the principal solution associated to the $\mu$-Beltrami equation is a bi-Lipschitz quasiconformal map.
}
\vspace{2mm}

We first show that it suffices to prove Theorem 2 for $N=1$, that is, a domain with one corner.  Suppose $\Omega$ is as assumed in Theorem 2 with corners of opening $\alpha_j$ at $z_j$ for $j=1,\ldots, N$.  Let $r_j>0$ be small enough so that $D(z_j,r_j)$ contains only one corner, $z_j$, and that $D(z_j,r_j)\cap \Omega$ connected.  Let $\Omega_j$ be a domain contained in $D(z_j,r_j)$ with one corner which coincides with $\Omega$ in a neighborhood of $z_j$.  Let $\Omega'$ be a $C^{1+\varepsilon}$ domain which coincides with $\Omega$ outside $\cup D(z_j,r_j)$ and such that the symmetric difference of $\Omega'$ and $\cup \Omega_j$ consists of $N+1$ connected components of positive pairwise distance.  Then  Theorem 1 and the $N=1$ version of Theorem 2 apply to   $\Omega'$ and $\Omega_j$ respectively, and we use the discussion in Section 2.2 to extend  these local results to the whole domain $\Omega$.


We now prove the $N=1$ version of Theorem 2 in two steps.  We first prove it for a special type of domain with one corner, and then use this version and conformal mapping to show the general case.

\subsection{Proof of Theorem 2 for Ice Cream Cone Domains}

Let $\L$ be a bounded domain in $\mathbb{C}$ which coincides with $\mathbb{S}_0^{\alpha \pi}$ for some $\alpha \in (0,2)$ in a neighborhood of 0, and whose boundary is $C^{1+\varepsilon}$ except at 0.  When $\alpha \in (0,1)$ one may think of the domain as an ice cream cone with a smoothly meeting smooth scoop of ice cream on top.  Hence we refer to these $\L$ as ice cream cone domains.  Let $\mu$ be a dilatation function supported on $\L$ with $\mu(0)$ allowable for $\mathbb{S}_0^{\alpha \pi}$ and $\varepsilon$-H\"older continuous.  We wish to show that $\Phi$, the quasiconformal map with dilatation $\mu$ with $\Phi(z) =z + O(1/z)$ near $\infty$, is bi-Lipschitz.

By Theorem 1 and the discussion in Section 2.2, it suffices to show that $\Phi$ is bi-Lipschitz in a neighborhood of 0, call such a neighborhood $U$.  Let $f_\angle$ be the bi-Lipschitz solution to $\mu_{f_\angle} = \mu(0) \cdot \chi_{\mathbb{S}_0^{\alpha \pi}}$ constructed in Section 2.5.  Then $f_\angle(U)$ takes $\mathbb{S}_0^{\alpha \pi}$ to a logarithmic spiral.  Consider
\[ h := \Phi \circ f_\angle^{-1}. \]
If we show that $h$ is bi-Lipschitz in $f_\angle(U)$, then $\Phi$ will be bi-Lipschitz in $U$.

By the composition formula for dilatations,   as $\mu$ is $\varepsilon$-H\"older continuous in $\L$,
$\mu_h$ is $\varepsilon$-H\"older continuous in $f_\angle(U)$. In particular
\[ |\mu_h (z)| \leq C |z|^\varepsilon \]
in $f_\angle(U)$  since $\mu_h (0) =0$.  Thus by the Lehto inequality (\ref{lehto}) for all annuli centered at 0 in $f_\angle(U)$,
\[ | M(A) - M(h(A)) | \leq C_1 \]
where $C_1$ does not depend on the specific annulus.  Using this inequality with $A = A(0,r,1)$ for small $r$, via the inequalities in (\ref{modineq}), we find some $M>1$ for which
\begin{equation}
 \frac1M |z| \leq |h(z)| \leq M |z| \label{0bilip}
\end{equation}
for $z \in f_\angle(U)$. We now let $z_1,z_2 \in f_\angle(U)$ be given and set $\epsilon < 1/M^2$.  If $|z_1| \leq \epsilon |z_2|$ or $|z_2| \leq \epsilon |z_1|$, we may use (\ref{0bilip}) to find
\[ \frac1{\tilde M} |z_1-z_2| \leq |h(z_1) - h(z_2)| \leq \tilde M |z_1 -z_2| \]
for $\tilde M = M^2/(M-1)$.

Now assume that $|z_1| > \epsilon |z_2|$ and $|z_2|>\epsilon |z_1|$.  We will make use of the self-similarity of the logarithmic spiral with respect to dilations.

Consider $f_\angle(\mathbb{S}_0^{\alpha \pi}) \cap A(0,\epsilon/2,1+2\epsilon)$.  There is a domain $\Delta$ with $C^{1+\varepsilon}$ boundary such that
\[ f_\angle(\mathbb{S}_0^{\alpha \pi}) \cap A(0,\epsilon/2,1+2\epsilon) \subset \Delta \subset f_\angle(\mathbb{S}_0^{\alpha \pi}). \]
As $\mu_h$ is an $\varepsilon$-H\"older continuous dilatation function supported on $f_\angle(\mathbb{S}_0^{\alpha \pi})$,
\[ \mu_\Delta = \mu_h \cdot \chi_\Delta \]
is an $\varepsilon$-H\"older continuous dilatation supported on a domain $\Delta$ with a $C^{1+\varepsilon}$ smooth  boundary.  So the principal solution corresponding to $\mu_\Delta$,  $\Phi_\Delta$, is bi-Lipschitz via Theorem 1. The bi-Lipschitz constant depends on $\| \mu \|_\infty$, $\mu$'s $\varepsilon$-H\"older constant, and the geometry of $\partial \Delta$.

Writing $z_1 = s e^{i \theta_1}$ and $z_2 = t e^{i \theta_2}$ with $t<s<1$ consider the map
\[ H(w) = \frac{h(sw)}s \]
for $w \in D(0, 1+2\epsilon)$.  Let 
\[ \tilde \mu_\Delta = \mu_H \chi_\Delta .\]
Then $\tilde \Phi_\Delta$, the principal $\tilde \mu_\Delta$-Beltrami solution, is bi-Lipschitz with constant dependent only on $\| \mu \|_\infty$, $\mu$'s $\varepsilon$-H\"older constant, and $ \partial \Delta$ as above with $\Phi_\Delta$.  As $A(0,\epsilon, 1+\epsilon)$ is compactly contained in $A(0,\epsilon/2, 1+2\epsilon)$, we get that $H$ is similarly bi-Lipschitz on $A(0,\epsilon, 1+\epsilon)$ to $\tilde \Phi_\Delta$ via Section 2.2.  So there exists an $N>1$ not dependent on $t,s,\theta_1$, or $\theta_2$ such that
\[ \frac{1}{N} |e^{i \theta_1} - \frac{t}{s} e^{i\theta_2}| \leq |H(e^{i \theta_1}) - H(\frac{t}{s} e^{i\theta_2})| \leq N |e^{i \theta_1} - \frac{t}{s} e^{i\theta_2}| \]
and multiplying by $s$ gives
\[ \frac1N |z_1 - z_2| \leq |h(z_1) - h(z_2)| \leq N |z_1-z_2| \]
for $z_1,z_2 \in f(U)$.  By choosing $L = \max(\tilde M,N)$ we have $h$ is $L$-bi-Lipschitz in $f_\angle(U)$.  Hence $\Phi = h \circ f_\angle^{-1}$ is bi-Lipschitz in $U$ as desired.  So Theorem 2 is valid for the ice cream cone domain $\L$.

\subsection{Proof of Theorem 2} 
Let $\Omega$ be a domain with a single $\varepsilon$-H\"older continuous corner of opening $\alpha \in (0,2)$ at, without loss of generality, 0 and $C^{1+\varepsilon}$ boundary elsewhere.  Let $\L$ be an ice cream cone domain with corner at 0 with opening $\alpha$.  Suppose we can construct a quasiconformal map $\Psi$, bi-Lipschitz in $\mathbb{C}$, taking $\L$ to $\Omega$, $\Psi(0) = 0$, with $\varepsilon$-H\"older continuous dilatation $\mu$ supported on $\L$ with $\mu(0)=0$.  Then we have reduced Theorem 2 to the case of ice cream cone domains in the same manner that we reduced Theorem 1 to the case of the unit disc.  We construct such a $\Psi$.

Let $\phi_\Omega$ ($\phi_\L$) be a conformal map from $\overline{\mathbb{D}}^c$ to $\overline{\Omega}^c$ $\left( \overline{\L}^c \right)$ which is $z+O(1)$ near $\infty$. By Theorem B, we have extensions for $\phi_\Omega$ and $\phi_\L$ so that
\[ \phi_\Omega \circ \phi_\L^{-1} : \L^c \to \Omega^c \] 
is $C^{1+\varepsilon}$ on the boundary and so is bi-Lipschitz on $\L^c$. We may assume each maps 1 to 0, and hence the singularities of the derivatives at the corner are canceled due to (\ref{cornerexp}).

Now let $\psi_\Omega$ ($\psi_\L$) be a conformal map from $\overline{\mathbb{D}}$ to $\overline{\Omega}$ ($\overline{\L}$) with 1 mapped to 0.  Again, by Theorem B,
\begin{equation} \psi_\Omega^{-1} \circ \phi_\Omega \circ \phi_\L^{-1} \circ \psi_\L
\label{comp}
\end{equation}
is a $C^{1+\varepsilon}$ homeomorphism of $\partial \mathbb{D}$ to itself.  We extend this homeomorphism precisely as we did in Section 3.1 above.  Call this extension $f$.  By freedom of choice of the particular conformal map $\psi_\Omega$, we may force $\mu_f (1) =0$ by making the derivative of (\ref{comp}) at 1 equal to 1. Set 
\[ \Psi := \left\{ \begin{array}{cr}
\phi_\Omega \circ \phi_\L^{-1} & \mbox{outside } \L\\
\psi_\Omega \circ f \circ \psi_\L^{-1} & \mbox{inside } \L 
\end{array} \right. .\]
This $\Psi$ is our desired map.

\end{document}